# ASYMPTOTIC EQUIVALENCE OF EMPIRICAL LIKELIHOOD AND BAYESIAN MAP

By Marian Grendár[1] and George Judge

*Bel University and University of California, Berkeley*

In this paper we are interested in empirical likelihood (EL) as a method of estimation, and we address the following two problems: (1) selecting among various empirical discrepancies in an EL framework and (2) demonstrating that EL has a well-defined probabilistic interpretation that would justify its use in a Bayesian context. Using the large deviations approach, a Bayesian law of large numbers is developed that implies that EL and the Bayesian maximum a posteriori probability (MAP) estimators are consistent under misspecification and that EL can be viewed as an asymptotic form of MAP. Estimators based on other empirical discrepancies are, in general, inconsistent under misspecification.

**1. Introduction.** Owen's empirical likelihood (EL) theorem ([30] and [31]) provides under traditional assumptions a basis for forming confidence regions for multivariate means and parameters in estimating equations. The basic EL idea is to proceed as if the sample $X_1, X_2, \ldots, X_n$, drawn from an unknown distribution $r(x)$ is i.i.d. and can be modeled as a multinomial distribution based on the observations. Inference for the $J$ unknown parameters is based on $K$ estimating equations and a nonparametric likelihood ratio statistic that asymptotically has a chi-square distribution. As a result, EL is an attractive orthodox semiparametric method of estimation and inference whose scope has been extended in several productive directions (e.g., see [17] and [41]).

Building on Owen's EL insight, this paper is concerned with using empirical likelihood as a method of estimation (cf. [3, 19, 27, 31] and [32], among others). Through a Bayesian law of large numbers (BLLN, Theorem

Received February 2008; revised June 2008.
[1]Supported by VEGA 1/3016/06 Grant.

2.2, Section 2.2, see also Section 2.3) we establish in a Bayesian setting an asymptotic connection between the EL and maximum a posteriori probability (MAP) estimators. The BLLN implies that under certain conditions the posterior measure asymptotically weakly concentrates on the MAP, or equivalently EL, estimators, even when the model is not correctly specified. The BLLN result is established via the large deviations (LD) approach (cf. Section 2.1).

Typically, as noted above, EL estimators are formed in an Estimating Equations (EE) framework [14]. In this way EL combines the flexibility of the nonparametric approach with the advantages of a finite parametrization. To be more specific, let us assume that a researcher is willing to specify only some of features of the data-sampling distribution $r(x;\theta)$. These model features, that is, the model $\Phi(\Theta)$, can be characterized by estimating equations (cf. [31], Chapter 3.5 and [27], Chapter 11): $\Phi(\Theta) \triangleq \bigcup_\Theta \Phi(\theta)$, where $\Phi(\theta) \triangleq \{q(x;\theta) : \int q(x;\theta) u_j(x;\theta) = 0, 1 \leq j \leq J; \int q(x;\cdot) = 1, q(x;\cdot) \geq 0\}$, $\theta \in \Theta \subseteq \mathbb{R}^K$. The $J$ number of estimating functions $u(\cdot)$ need not be equal to the number $K$ of parameters. Given the sample drawn from $r(x;\theta)$, the EL estimator $\hat\theta_{\mathrm{EL}}$ of $\theta$ is obtained as a parametric component of

$$(1.1) \qquad \hat q_{\mathrm{EL}}(x_i; \hat\theta_{\mathrm{EL}}) = \arg \sup_{q(x_i;\theta) \in \widehat\Phi(\theta)} \sup_{\theta \in \Theta} \sum_{i=1}^n \log q(x_i;\theta),$$

where $\widehat\Phi(\theta) \triangleq \{q(x_i;\theta) : \sum_{i=1}^n q(x_i;\theta) u_j(x_i;\theta) = 0, 1 \leq j \leq J; \sum_{i=1}^n q(x_i;\theta) = 1, q(x_i;\cdot) \geq 0\}$ is the empirical form of $\Phi(\theta)$.

The asymptotic performance of $\theta_{\mathrm{EL}}$ has been studied by Qin and Lawless [32]. The same asymptotic properties are exhibited by the exponential tilting estimator, which results when in (1.1) the log-likelihood is replaced by the negative Kullback–Leibler discrepancy (empirical entropy) $-\sum_{i=1}^n q(x_i;\theta) \log(q(x_i;\theta)/\mu(x_i))$ of $q$ with respect to the uniform probability mass function $\mu(\cdot)$ (cf. [18, 23] and [27]). Also, in [3], the Euclidean distance $\sum_{i=1}^n (q(x_i;\theta) - \mu(x_i))^2$ is employed and in recent years the Cressie–Read [4] family of discrepancy measures has been used in the EL framework. When the model is correctly specified, the resulting estimators are consistent regardless of the discrepancy measure used. Since in practice the model is rarely specified correctly, it is of interest to study consistency of various EL estimators under statistical model misspecification.

Recognizing the importance of the statistical model in the estimation process, we study consistency under misspecification in the Bayesian setting. By means of the LD approach we obtain the Bayesian law of large numbers (Theorem 2.2). The BLLN together with Lemmas 2.1 and 2.2 imply that, in the Bayesian setting, MAP and EL estimators are consistent under misspecification. In the Bayesian setting, Euclidean and other nonlikelihood members of the Cressie–Read family may, in general, be inconsistent. The same holds for the posterior mean (cf. Section 2.3, Example 2.1).

The BLLN also sheds a new light on the problem of extending of EL into a Bayesian method (cf. [26, 31, 33] and [35]). In [26] Lazar listed possible ways of turning EL into a Bayesian method and studied one of the possibilities within the framework of Monahan and Boos [28]. Schennach [35] proposed a specific prior over a set of sampling distributions to get a Bayesian procedure that admits an operational form *similar* to EL. In [33] a different prior over the set of probability measures is considered and a group of EL-like methods is obtained. The BLLN together with Lemmas 2.1 and 2.2 imply that if certain requirements are put on an infinite-dimensional prior, the Bayesian MAP method [cf. equation (2.2)] leads asymptotically to the same point estimator(s) as empirical likelihood. Informally put, this means that EL, as a method of estimation, can be viewed as an asymptotic instance of MAP. Extension of the connection between EL and MAP into the field of inference remains an open problem; compare, for instance, Freedman's [10], where it is shown that the distributional asymptotic properties of maximum nonparametric likelihood (a special, nonparametric case of EL) and MAP estimators might be different, since the Bernstein–von Mises theorem does not apply, even for a simple infinite dimensional models.

1.1. *Organization of the paper.* In Section 2.1 the LD approach to Bayesian consistency is informally described. In Section 2.2 a basic framework is established, Bayesian nonparametric consistency is formally defined, $L$-divergence is introduced and the BLLN theorem is proved for the i.i.d. case. In Section 2.3 the BLLN for the semiparametric model is discussed and it is demonstrated that the $L$-projection, singled out by the BLLN, is an asymptotic form of the EL and MAP estimators. In order to further explore consistency under misspecification and expand the scope of the related asymptotic connection between MNPL/EL and MAP, we prove the BLLN also for the multicolor Pólya sampling process (Section 2.4) where using the BLLN suggests two possible variants of MNPL. In Section 2.5, the BLLN is proved for right censored data, showing that the Kaplan–Meier estimator is an asymptotic form of Bayesian MAP.

**2. Bayesian LLN's, MNPL, EL and MAP.** In general, a feasible set $\Phi(\Theta)$ of nonparametric sampling distributions which are indexed by a parameter $\theta$ can be formed in a way different from the conventional EE described above. The purely nonparametric $\Phi$ is contained in $\Phi(\Theta)$ as a special case.

2.1. *Large deviations approach to Bayesian consistency under misspecification.* In a Bayesian framework, a prior distribution $\Pi$ over the set $\Phi(\Theta)$ is assumed, and it induces a prior distribution $\Pi(\theta)$ over $\Theta$. Assuming the Bayesian framework, it is of interest to know the sampling distribution(s) on which the posterior measure concentrates, as $n$, the sample size gets large. The importance of the frequentist concept of consistency for Bayesian statistics can be justified both from subjectivist and objectivist Bayesian positions (cf. [40] or [13], Chapter 4). A formal definition of Bayesian consistency is given in Section 2.2.

We use a large deviations (LD) (cf. Ben-Tal, Brown and Smith [1] and [2] and Ganesh and O'Connell [11]) approach to Bayesian nonparametric consistency. The approach results in a Bayesian–Sanov Theorem (BST) and its corollary the Bayesian law of large numbers (BLLN) that establishes the consistency. LD theory is a subfield of probability theory where, informally, the typical concern is about the asymptotic behavior, on a logarithmic scale, of the probability of a given event. The BST identifies the rate function governing exponential decay of the posterior measure, and this in turn identifies the sampling distributions on which the posterior concentrates, as those distributions that minimize the rate function. Currently used approaches to Bayesian nonparametric consistency (cf. [38]) do not recognize this concentration of the posterior measure as a solution of the optimization problem.

The Bayesian law of large numbers may be, informally, stated as follows: if the prior over a set $\Phi$ of sampling distributions, which might not include the "true" distribution with probability density function $r$, satisfies certain conditions, then the posterior asymptotically concentrates (a.s. $r^\infty$) on weak neighborhoods of the $L$-projections of $r$ on $\Phi$. $L$-projection $\hat{q}$ of $r$ on $\Phi$ is $\hat{q} = \arg\inf_{q \in \Phi} L(q||r)$, where $L(q||r)$ is the $L$-divergence of probability density function $q$ with regard to $r$. In the case of i.i.d. sampling, $L(q||r) = -\int r \log q$.

Finally, let us note that the BST is Bayesian counterpart of a Sanov theorem for empirical measures (cf. [34] and [5], Sections III and VII and references cited therein). The latter, as well as its corollary, the conditional law of large numbers, are basic results of large deviations (LD) theory (cf. [5] and [7]). The LD theorems for empirical measures have a bearing for the relative entropy maximization method. Kitamura and Stutzer [22] noted that the LD argument can be used also in the semiparametric EE setting, where it provides an underpinning to exponential tilting (see [18] and [22]), also known as the maximum entropy empirical likelihood [27] method. In fact, the work of Kitamura and Stutzer [23] served as a starting point for our attempt to provide a similar underpinning to the MNPL and EL methods (see also [16]). It turned out that this is only possible in a Bayesian framework.

2.2. *BLLN for i.i.d. sampling.* Let $\mathcal{P}$ be the set of all probability measures on $(\mathbb{R}, \mathcal{B})$, which are dominated by the Lebesgue measure. Let $X_1, X_2, \ldots$ be i.i.d. random variables that take values in $(\mathbb{R}, \mathcal{B})$, with probability density function (PDF) $r$ where probability densities are denoted by lower case. $\mathcal{P}$

is endowed with weak topology. Let $\Phi \subseteq \mathcal{P}$. It is not assumed that $r$, the true sampling distribution, is necessarily in $\Phi$. Let $\sigma(\mathcal{P})$ be a Borel $\sigma$-field on $\mathcal{P}$. A positive prior $\Pi$ is put on $(\mathcal{P}, \sigma(\mathcal{P}))$ that is strictly positive over $\Phi$. The prior combines with data $X_1^n \triangleq X_1, X_2, \ldots, X_n$ to define the posterior distribution

$$\Pi_n(Q|X_1^n) = \frac{\int_Q e^{-l_n(q)} \Pi(dq)}{\int_\Phi e^{-l_n(q)} \Pi(dq)},$$

where $l_n(q) \triangleq -\sum_{i=1}^n \log q(x_i)$, $Q \subseteq \Phi$.

Let $d$ be a metric on $\mathcal{P}$. The sequence $\{\Pi_n(\cdot|X_1^n), n \geq 1\}$ is said to be $d$-consistent at $r$, if there exists a $\Omega_0 \subset \mathbb{R}^\infty$ with $r(\Omega_0) = 1$ such that for $\omega \in \Omega_0$, for every neighborhood $U$ of $r$, $\Pi_n(U|X_1^n) \to 1$ as $n$ goes to infinity. If a posterior is $d$-consistent for any $r \in \Phi$, then it is said to be $d$-consistent. If the consistency holds for the Hellinger distance, then the posterior is strongly consistent. If convergence holds in weak topology, the posterior is said to be weakly consistent. In [39] a decision-theoretic argument is proposed in favor of weak consistency. Surveys of Bayesian nonparametric consistency can be found in [12, 13] and [40].

To the best of our knowledge, Ben-Tal, Brown and Smith [1] were the first to use an LD approach to Bayesian nonparametric consistency. The authors showed consistency for $X$, taking values from a finite set $\mathcal{X}$ and under a possibly misspecified model. Recently, Ganesh and O'Connell [11] independently established the first formal BST, for finite set $\mathcal{X}$ and a well-specified model. Here we develop the BST and the BLLN for $\mathcal{X} = \mathbb{R}$ and a possibly misspecified model. Using techniques other than LD, consistency in the Hellinger distance and under misspecification was studied by Kleijn and van der Vaart [24].

The key quantity that governs the LD exponential decay of the posterior $\Pi_n(Q|X_1^n)$ in the i.i.d. case is the $L$-divergence of $q \in \mathcal{P}$ with regard to $p \in \mathcal{P}$: $L(q\,||\,p) \triangleq -\int p \log q$ (cf. [15]). In the discrete case, $L$-divergence appears in Freedman's ([9], Theorem 1) as "entropy." If $p$ is an empirical PMF, then $L$-divergence appears as Kerridge's inaccuracy ([21] and [25]) which is just the negative of the nonparametric likelihood. The $L$-projection $\hat{q}$ of $p$ on $Q \subseteq \mathcal{P}$ is $\hat{q} \triangleq \arg\inf_{q \in Q} L(q||p)$. The value of $L$-divergence, at an $L$-projection of $p$ on $Q$, is denoted by $L(Q||p)$.

Finally, let, for $p, q \in \mathcal{P}$, $\varepsilon > 0$, $B_\varepsilon(q, p) \triangleq \{q' \in \mathcal{P} : L(q'\,||\,p) - L(q\,||\,p) < \varepsilon\}$. For $A \subseteq \mathcal{P}$, $B_\varepsilon(A, p) \triangleq \{q \in \mathcal{P} : L(q\,||\,p) - L(A\,||\,p) < \varepsilon\}$.

Using this notation, the BST can be stated as follows:

THEOREM 2.1 (BST). *Let $X_1^n$ be i.i.d. $r$. Let $Q$ and $\Phi$ be open in weak topology; $Q \subset \Phi \subseteq \mathcal{P}$. Let $L(Q\,||\,r) < \infty$; for any $\varepsilon > 0$, let $\Pi(B_\varepsilon(Q, r)) > 0$ and $\Pi(B_\varepsilon(\Phi, r)) > 0$. Then, for $n \to \infty$,*

$$\frac{1}{n} \log \Pi_n(q \in Q|X_1^n) = -\{L(Q\,||\,r) - L(\Phi\,||\,r)\} \qquad a.s. \ r^\infty.$$

PROOF. For $S \subseteq \mathcal{P}$, $l_n(S) \triangleq \inf_{q \in S} l_n(q)$. Let for $\varepsilon > 0$, $B_\varepsilon^n(S) \triangleq \{q : l_n(q) - l_n(S) < \varepsilon\}$. Then, $\int_A e^{-l_n(A)} \Pi(dq)$, $A = \{Q, \Phi\}$, can be bounded as

$$e^{-l_n(A)-\varepsilon} \Pi(B_\varepsilon^n(A) \cap A) \leq \int_A e^{-l_n(q)} \Pi(dq) \leq e^{-l_n(A)}.$$

By lower semicontinuity of $L$-divergence in the weak topology and a strong law of large numbers [which can be applied, since $L(Q\,||\,r) < \infty$, by assumption], $\frac{1}{n} l_n(A) \to L(A\,||\,r)$, a.s. $r^\infty$, as $n \to \infty$. Thus, it holds:

$$\limsup_{n \to \infty} \frac{1}{n} \log \int_A e^{l_n(q)} \Pi(dq) \leq -L(A\,||\,r).$$

So, $\limsup_{n \to \infty} \Pi_n(Q|X_1^n) \leq -\{L(Q\,||\,r) - L(\Phi\,||\,r)\}$. By the same argument (SLLN and continuity), for sufficiently large $n$, $\Pi(B_\varepsilon^n(A)) > 0$, since $\Pi(B_\varepsilon(A, r)) > 0$ by assumption. As $B_\varepsilon^n(A) \cap A \neq \varnothing$, thus $\lim_{n \to \infty} \frac{1}{n} \log \Pi(B_\varepsilon^n(A) \cap A) = 0$. Hence, $\liminf_{n \to \infty} \Pi_n(Q|X_1^n) \geq -\{L(Q\,||\,r) - L(\Phi\,||\,r)\}$. □

The posterior probability $\Pi_n(Q|X_1^n)$ decays exponentially fast with the decay rate $L(Q\,||\,r) - L(\Phi\,||\,r)$. The BST implies the Bayesian law of large numbers (BLLN).

THEOREM 2.2 (BLLN). *Let $\Phi \subseteq \mathcal{P}$ be open in weak topology. Let (1) for every $q \in \Phi$, $\Pi(B_\varepsilon(q,r)) > 0$ and (2) $L(\Phi\,||\,r) < \infty$. Let $U \triangleq \bigcup_k W(\hat{q}_k, \varepsilon)$ be a union of weak $\varepsilon$-balls $W(\hat{q}_k, \varepsilon)$ centered at $L$-projections $\hat{q}_k$, $k = 1, \ldots, \kappa$, $\kappa < \infty$, of $r$ on $\Phi$. Then,*

$$\lim_{n \to \infty} \Pi_n(q \in U|X_1^n) = 1 \qquad a.s.\ r^\infty.$$

PROOF. Let $Q \subset \Phi$ be open sets in weak topology. First, let $Q$ be any set such that $\infty > L(Q\,||\,r) > L(\Phi\,||\,r)$. Then, assumptions of the BST are satisfied, and the theorem implies that $\Pi_n(Q|X_1^n) \to 0$, a.s. $r^\infty$, as $n \to \infty$. Note that $L(Q\,||\,r) = \infty$ for such $Q$ that the $L$-projection $\hat{q}^Q$ of $r$ on $Q$ has support that is smaller than the support of $r$. However, for such a $\hat{q}^Q$, the posterior probability would be zero. The posterior thus concentrates on $L$-projections of $r$ on $\Phi$, provided that their support is not smaller than that of $r$. This is guaranteed by the assumption $L(\Phi\,||\,r) < \infty$. □

The BLLN theorem is an extension of Schwartz' consistency theorem [36], to the case of a misspecified model. Assumptions 1 and 2 of the BLLN, called hereafter the Schwartz conditions, reduce in the well-specified case to the Kullback–Leibler support condition (cf. [36] and [13], Theorem 4.4.2).

The next lemma points out that the Bayesian maximum a posteriori probability (MAP), which selects $\hat{q}_{\text{MAP}} \triangleq \arg\sup_{q \in \Phi} \Pi_n(q\,|\,X_1^n)$, satisfies the BLLN.

LEMMA 2.1. *Let $\Phi \subseteq \mathcal{P}$ be open and the Schwartz conditions be satisfied. Then, as $n \to \infty$, the set of MAP distributions $\mathcal{M} \triangleq q\{\hat{q}_{\mathrm{MAP}} : \hat{q}_{\mathrm{MAP}} = \arg\sup_{q \in \Phi} \Pi_n(q \,|\, X_1^n)\}$ converges (a.s. $r^\infty$) to the set of L-projections of $r$ on $\Phi$.*

PROOF. Thanks to the Strong LLN (SLLN), which can be applied under the Schwartz condition 2, the conditions for the infimum of minus the logarithm of the posterior probability (positivity of which is guaranteed by the Schwartz condition 1) turn into those for $L$-projections. □

Directly from the strong LLN it follows that the maximum nonparametric likelihood (MNPL), that selects $\hat{q}_{\mathrm{MNPL}} \triangleq \arg\inf_{q \in \Phi} l_n(q)$, satisfies the BLLN.

LEMMA 2.2. *Let $\Phi \subseteq \mathcal{P}$ be open and Schwartz condition 2 be satisfied. Then, as $n \to \infty$, the set of MNPL distributions converges (a.s. $r^\infty$) to the set of L-projections of $r$ on $\Phi$.*

Selection of a posterior mean or a sampling distribution that minimizes, say, the Kullback–Leibler distance $I(q \,||\, r) = \int q \log \frac{q}{r}$ with regard to $q$, in a misspecified case, would in general violate the BLLN.

The lemmas also mean that the MNPL and the MAP methods asymptotically select the same sampling distribution(s).

Next, we turn to the semiparametric setting.

2.3. *BLLN for the semiparametric $\Phi(\Theta)$.* Let $X$ be a random variable with probability density function $r(X; \theta)$ parametrized by $\theta \in \Theta \subseteq \mathbb{R}^K$. A Bayesian specifies a model $\Phi(\Theta)$ and puts a positive prior $\Pi$ over $\Phi(\Theta)$, which in turn induces a prior $\Pi(\theta)$ over $\Theta$; see Florens and Rolin [8], where also several models are worked out using a Dirichlet process prior. If the requirements of the BLLN are satisfied, then the posterior $\Pi_n(\cdot|X_1^n)$ concentrates on weak neighborhoods of $L$-projections $\hat{q}$ of $r$ on $\Phi(\Theta)$,

$$\hat{q}(x; \hat{\theta}) = \arg \inf_{q(x;\theta) \in \Phi(\theta)} \inf_{\theta \in \Theta} L(q(x;\theta) \,||\, r).$$

The most common form of $\Phi(\Theta)$ is the one defined by estimating equations (cf. Section 1). In this case, $\Phi(\Theta)$ is also known as a linear family of distributions that we denote as $\mathcal{L}(u, \theta)$. The $L$-projection of $r$ on $\mathcal{L}(u, \theta)$ can be found by means of the following Theorem 2.3. To state it, we introduce a $\Lambda$ family of distributions and recall the concept of support of a convex set. Let $\Lambda$ be a family of probability density functions: $\Lambda(r, u, \lambda, \theta) \triangleq \{p \in \mathcal{P} : p = r[1 - \sum_{j=1}^J \lambda_j u_j(\cdot; \theta)]^{-1}, \lambda \in \mathbb{R}^J\}$. The support $S(\mathcal{C})$ of a convex set $\mathcal{C} \subset \mathcal{P}$ is just the support of the member of $\mathcal{C}$ for which $S(\cdot)$ contains the support of any other member of the set.

THEOREM 2.3. *Let $\Phi = \mathcal{L}(u, \theta)$. Let $r \in \mathcal{P}$ be such that $S(r) = S(\mathcal{L})$. Then, the L-projection $\hat{q}$ of $r$ on $\Phi$ is unique and belongs to the $\Lambda(r, u, \lambda, \theta)$ family; that is, $\mathcal{L}(u, \theta) \cap \Lambda(r, u, \lambda, \theta) = \{\hat{q}\}$.*

PROOF. In light of Theorem 9 of [6] it suffices to check that $\hat{q} = r[1 - \sum_{j=1}^{J} \lambda_j u_j(\cdot; \theta)]^{-1}$, with $\lambda$ such that $\hat{q} \in \mathcal{L}(u, \theta)$, satisfies $\int_{S(r)} r(1 - \frac{q'}{\hat{q}}) = 0$, for all $q' \in \Phi$, which is indeed the case. □

The estimator $\hat{\theta}$ can, thanks to convex duality, be obtained as

$$\hat{\theta} = \arg \inf_{\theta \in \Theta} \sup_{\lambda \in \mathbb{R}^J} L(q(x; r, u, \lambda, \theta) \| r),$$

where $q(x; r, u, \lambda, \theta) \in \Lambda(r, u, \lambda, \theta)$. Since $r$ is in practice not known, Kitamura and Stutzer [22] suggested that $L(q(x; r, u, \lambda, \theta) \| r)$ be replaced by its estimate $\hat{L}(q(x; \mu, u, \lambda, \theta)) \triangleq -\sum_{i=1}^{n} \log q(x_i; \mu, u, \lambda, \theta)$, where $q(x; \mu, u, \lambda, \theta)$ belongs to $\Lambda(\mu, u, \lambda, \theta)$ and $\mu$ is the uniform PMF on $X_1^n$. The resulting estimator

(2.1) $$\hat{\theta}_{\mathrm{EL}} \triangleq \arg \inf_{\theta \in \Theta} \sup_{\lambda \in \mathbb{R}^J} \hat{L}(q(x; \mu, u, \lambda, \theta))$$

is just the empirical likelihood (EL) estimator ([32] and [31]), since (2.1) is a convex dual problem to the optimization problem (1.1), by means of which EL is usually defined. Analogously to Lemma 2.2, it can be shown that the EL estimator $\hat{q}_{\mathrm{EL}}(x; \hat{\theta}_{\mathrm{EL}})$ asymptotically (a.s. $r^\infty$) turns into an L-projection of $r$ on $\Phi(\Theta)$. The same holds for the MAP estimator

(2.2) $$\hat{q}_{\mathrm{MAP}}(x; \hat{\theta}_{\mathrm{MAP}}) = \arg \sup_{q(x; \theta) \in \Phi(\Theta)} \sup_{\theta \in \Theta} \Pi_n(q(x; \theta) | x_1^n).$$

Hence, the EL and the MAP estimators are consistent under misspecification. This provides a basis for the EL approach as well for the Bayesian MAP estimation.

EL estimators which are based on other discrepancy measures are, in general, not consistent when the model is not correctly specified. Example 2.1 illustrates the inconsistency of the posterior mean.

EXAMPLE 2.1. Let $\Phi = \Phi_1 \cup \Phi_2$, where $\Phi_1 = \mathcal{L}(u, \theta)$ and $\theta \in \Theta_1 = \{\theta \in \mathbb{R} : \theta \leq \theta_1\}$. Similarly, $\Phi_2 = \mathcal{L}(u, \theta)$ and $\theta \in \Theta_2 = \{\theta \in \mathbb{R} : \theta \geq \theta_2\}$. Let $u(x, \theta) = x - \theta \in \mathbb{R}$. For a sampling distribution $r$, it is possible to find $\theta_2 > \mathrm{E}X > \theta_1$, such that $L(\Phi_1 \| r) = L(\Phi_2 \| r)$. Assume this to be the case. Then, under the Schwartz conditions, the posterior concentrates on the weak balls centered at the L-projections of $r$ on $\Phi_1$ and $\Phi_2$, rendering the posterior mean inconsistent.

In the univariate case, BST and BLLN (Theorems 2.1 and 2.2), $X$ can be replaced by a multivariate random variable and the theorems remain valid. Consequently, the extension to $\Phi$, constructed by multivariate EE, is also direct. As an example, consider the linear statistical model $Y = \alpha + \beta X + \varepsilon$, with stochastic $X$. In EE this is usually approached through estimating equations $\Phi(\theta) = \{q(x, y; \theta) : \int q(x, y; \theta)[Y - (\alpha + \beta X)] = 0, \int q(x, y; \theta) X[Y - (\alpha + \beta X)] = 0\}$ and $\theta \triangleq (\alpha, \beta) \in \mathbb{R}^2 \equiv \Theta$, which are based on the Gaussian model score equations. The multivariate BLLN shows that the posterior asymptotically concentrates on the $L$-projections of $r$ on $\Phi = \bigcup_\Theta \Phi(\theta)$, and EL and MAP comply with the BLLN.

2.4. *BLLN for Pólya sampling.* In this section we prove the BST and the BLLN for a multi-color Pólya urn—a simple sampling process where data are neither identically nor independently distributed. The theorems can also be directly used in a corresponding semiparametric $\Phi(\Theta)$ setting.

The probability of a sample $X_1^n$ being drawn from a multicolor Pólya urn, with parameter $c \in \mathbb{Z}$ and initial configuration $q(N) \triangleq (\alpha_1, \ldots, \alpha_m)/N$, is $\log \Pi(X_1^n \,|\, q(N); c) \triangleq \sum_{i=1}^m \sum_{l=0}^{n_i - 1} [\log(\alpha_i + jc) - \log(N + jc)]$; this is meaningful if $-nc \leq \min(\alpha_1, \ldots, \alpha_m)$. We embed the sampling scheme into a Bayesian nonparametric setting. To this end, let $\mathcal{P}(\mathcal{X})$ be set of all PMFs with the support $\mathcal{X} = \{x_1, \ldots, x_m\}$. Let $\Phi \subseteq \mathcal{P}(\mathcal{X})$ and let $\Phi(N)$ denote the intersection of $\Phi$ with the set of all possible configurations of the $N$-urn. Let $\Phi(N)$ be the support of the prior distribution $\Pi(q(N))$ of initial configurations $q(N)$. Let $r(N)$ be the true initial configuration, where $r(N)$ is not necessarily in $\Phi(N)$. As before, we are interested in the LD asymptotics of the posterior distribution $\Pi_n(q(N) \,|\, X_1^n; c)$. Asymptotic investigations of posterior consistency will be carried on under the following assumptions: (1) $n$ and $N$ go to infinity in such a way that $\beta(n) \triangleq \frac{n}{N} \to \beta \in (0, 1)$ as $n \to \infty$, (2) and $r(N)$ converges in the total variation metric to $r \in \mathcal{P}(\mathcal{X})$ as $n \to \infty$. Topological qualifiers are meant in the topology induced on the $m$-dimensional simplex by the usual topology on $\mathbb{R}^m$.

The exponential decay of the posterior is governed by Pólya $L$-divergence. For $p, q \in \mathcal{P}(\mathcal{X})$, the Pólya $L$-divergence $L_\beta^c(q \,||\, p)$ of $q$ with respect to $p$ is

$$L_\beta^c(q \,||\, p) \triangleq -\sum_{i=1}^m p_i \log(q_i + \beta c p_i) + \frac{1}{\beta c} \sum_{i=1}^m q_i \log \frac{q_i}{q_i + \beta c p_i}.$$

By the continuity argument, $L_\beta^0(q \,||\, p) \triangleq -\sum_{i=1}^m p_i \log q_i - 1$. The Pólya $L_\beta^c$-projection $\hat{q}$ of $p$ on $Q \subseteq \mathcal{P}(\mathcal{X})$ is $\hat{q} \triangleq \arg\inf_{q \in Q} L_\beta^c(q \,||\, p)$. The value of $L_\beta^c$-divergence at an $L_\beta^c$-projection of $p$ on $Q$ is denoted by $L_\beta^c(Q || p)$.

THEOREM 2.4. *Let $Q \subset \Phi$ be an open set. Let $\beta(n) \to \beta \in (0,1)$ and $r(N) \to r$ as $n \to \infty$. Let $L_\beta^c(Q\|r) < \infty$. Then, for $n \to \infty$,*

$$\frac{1}{n}\log \Pi_n(q(N) \in Q \,|\, X_1^n; c) = -\{L_\beta^c(Q\,\|\,r) - L_\beta^c(\Phi\,\|\,r)\}$$

*with probability one.*

PROOF. The proof is constructed separately for $c > 0$, $c < 0$ and $c = 0$.

For $c \neq 0 \wedge \eta c \notin (\mathbb{Z}^-)^m \wedge \eta \notin \mathbb{Z}^-$, $\log \Pi(X_1^n \,|\, q(N); c)$ can equivalently be expressed as $\log(\Gamma(\eta)/\Gamma(\eta+n)) + \sum_{i=1}^m \log(\Gamma(\eta q_i + n_i)/\Gamma(\eta q_i))$, where $\Gamma(\cdot)$ is the Gamma function and $\eta \triangleq N/c$. For $0 < a < b$, the ratio $\Gamma(b)/\Gamma(a)$ can be upper-bounded by $b^{b-1/2}/a^{a-1/2}e^{b-a}$ and lower-bounded by $b^{b-1}/a^{a-1}e^{b-a}$ (cf. [20]). Then, $\Pi_n(q(N) \in Q \,|\, x^n; c)$ can be upper-bounded by $U_n$ (dependence of $q$ on $N$ is made implicit),

$$U_n = \frac{\sum_{q \in Q} \Pi(q) \prod_{i=1}^m e^{-n\,l(q_i, 1/(2n))}}{\sum_{q \in \Phi} \Pi(q) \prod_{i=1}^m e^{-n\,l(q_i, 1/n)}},$$

lower-bounded by $L_n$ in similar way; to get $L_n$ just replace $1/2n$ with $1/n$ in $U_n$. There, $l(q_i, \alpha) \triangleq -[(\gamma_i - \alpha)\log \gamma_i + (\gamma_i + \nu_i^n - \alpha)\log(\gamma_i + \nu_i^n)]$, $\gamma_i \triangleq \frac{q_i}{\beta(n)c}$, $\alpha \in \{\frac{1}{n}, \frac{1}{2n}\}$ and $\nu^n$ is the empirical measure induced by the sample $X_1^n$.

Next, we use simple bounds to upper bound $U_n$ by $\bar{U}_n$

$$\bar{U}_n = \frac{\prod_{i=1}^m e^{-n\,l(\hat{q}_i(Q, 1/(2n)), 1/(2n))}}{\pi(\hat{q}(\Phi, 1/n)) \prod_{i=1}^m e^{-n\,l(\hat{q}_i(\Phi, 1/n), 1/n)}},$$

and to lower bound $L_n$ by $\underline{L}_n$; to get $\underline{L}_n$ just replace $1/2n$ with $1/n$ in $\bar{U}_n$. There, $\hat{q}(\cdot, \alpha) \triangleq \arg\inf_{q \in \cdot} \sum_{i=1}^m l(q_i, \alpha)$.

By the Strong Law of Large Numbers for Pólya Sampling, $\nu^n \to r$, almost surely, as $n \to \infty$. The Pólya $L$-divergence is continuous in $q$ and $Q$ is open, by assumption. Thus, $\frac{1}{n}\log \bar{U}_n$ converges, with probability one, to $-\{L_\beta^c(Q\,\|\,r) - L_\beta^c(\Phi\,\|\,r)\}$, as $n \to \infty$. This is the same as the "point" of almost sure convergence of $\frac{1}{n}\log \underline{L}_n$ and the theorem for $c > 0$ is thus proven.

For $c \neq 0 \wedge (1-\eta q) \notin (\mathbb{Z}^-)^m \wedge (1-\eta) \notin \mathbb{Z}^-$, $\log \Pi(X_1^n \,|\, q(N); c)$ can equivalently be expressed as $\log(\Gamma(1-\eta-n)/\Gamma(1-\eta)) + \sum_{i=1}^m \log(\Gamma(1-\eta q_i)/\Gamma(1-\eta q_i - n_i))$. The proof then can be constructed along the same lines as for $c > 0$. The case of $c = 0$ is straightforward. □

From the Pólya BST (Theorem 2.5), the BLLN for Pólya sampling directly follows. It is worth noting that the MNPL in Pólya sampling can be constructed in two ways: either via maximization of $\Pi(X_1^n \,|\, q(N); c)$, or by maximization of the negative of $L_\beta^c(q\|\nu^n)$ with regard to $q$, where $\nu^n$ is empirical PMF induced by sample $X_1^n$. The methods could be called "exact" and "asymptotic" MNPL, respectively. Both the methods comply with the Pólya BLLN, as does the Bayesian MAP.

2.5. *BST for right-censored data.* Right-censoring of a r.v. $X$ by a r.v. $Y$ [both on $(\mathbb{R}, \mathcal{B})$] can be described by the following hierarchical model: $\delta \sim \text{Ber}(\alpha)$, $\alpha \triangleq \int F_0(y) \, dG_0(y)$; if $\delta = 0$, then $X \sim F_0$; if $\delta = 1$, then $X = (Y, \infty)$ where $Y \sim G_0$; $X$'s are conditionally independent. A Bayesian puts positive prior over the set $\Phi$ of distributions of $X$. Let the prior over distributions of $Y$ be concentrated at $G_0$. We are interested in the exponential decay of the posterior

$$\Pi_n(F \in Q \mid X_1^n) = \frac{\int_Q e^{-l_n(F, n_1)} \Pi(dF)}{\int_\Phi e^{-l_n(F, n_1)} \Pi(dF)},$$

where $l_n(F, n_1) \triangleq -\sum_{i:\delta_i=0} \log F(\{X_i\}) - \sum_{i:\delta_i=1} \log F((Y_i, \infty))$, $Q \subset \Phi$, $F_0$ is not necessarily in $\Phi$, and $n_1$ is the number of noncensored data, out of $n$ observations. The decay is governed by the $L$-divergence of $F$ with regard to $(F_0, G_0)$ for right-censoring

$$L(F \,||\, (F_0, G_0))$$
$$\triangleq - \left[ \alpha \int \log F(x) \, dF_0(x) + (1 - \alpha) \int \log F((y, \infty)) \, dG_0(y) \right].$$

The $L$-projection $\hat{F}$ of $(F_0, G_0)$ on $Q \subseteq \mathcal{P}$ is $\hat{F} \triangleq \arg\inf_{F \in Q} L(F \,||\, (F_0, G_0))$, and $L(Q \,||\, (F_0, G_0))$ denotes the value of the $L$-divergence at an $L$-projection of $F_0$ on $Q$. Let $B_\varepsilon(Q, F_0) \triangleq \{F \in \mathcal{P} : L(F \,||\, (F_0, G_0)) - L(Q \,||\, (F_0, G_0)) < \varepsilon\}$. The BST for right-censoring follows.

THEOREM 2.5. *Let $X_1^n$ be right-censored data generated by the above model. Let $Q, \Phi$ be open in weak topology; $Q \subset \Phi \subseteq \mathcal{P}$. Let $L(Q \,||\, (F_0, G_0)) < \infty$, and for any $\varepsilon > 0$, let $\Pi(B_\varepsilon(Q, F_0)) > 0$ and $\Pi(B_\varepsilon(\Phi, F_0)) > 0$. Then for $n \to \infty$,*

$$\frac{1}{n} \log \Pi_n(F \in Q \mid X_1^n) = -\{L(Q \,||\, (F_0, G_0)) - L(\Phi \,||\, (F_0, G_0))\}.$$

PROOF. Note that $\frac{1}{n} l_n(F, n_1)$ converges to $L(F \,||\, (F_0, G_0))$, with probability 1, by the SLLN. Arguments go along the lines of the proof of Theorem 2.1. □

From the BST (Theorem 2.5), the BLLN follows for right-censored data in the same way as it does for the i.i.d. case from Theorem 2.1. The BLLN for right-censoring demonstrates that the posterior concentrates on weak neighborhoods of the $L$-projections of $(F_0, G_0)$ on $\Phi$, if the $\varepsilon$-balls $B_\varepsilon(F, (F_0, G_0)) \triangleq \{F' \in \mathcal{P} : L(F' \,||\, (F_0, G_0)) - L(F \,||\, (F_0, G_0)) < \varepsilon\}$, have positive prior probability. This, together with assumption $L(\Phi \,||\, (F_0, G_0))$, forms the Schwartz conditions for right censoring.

Under the Schwartz conditions, a set of Bayesian MAP estimators $\hat{F}_{\text{MAP}} \triangleq \arg\sup_{F \in \Phi} \Pi_n(F \mid X_1^n)$ asymptotically coincides with a set of $L$-projections of $(F_0, G_0)$ on $\Phi$. The same holds true for the MNPL/EL estimator $\hat{F}_{\text{EL}} \triangleq \arg\inf_{F \in \Phi} l_n(F, n_1)$. The Kaplan–Meier estimator follows from $\hat{F}_{\text{EL}}$ in the standard way (cf. [31]). Thus, the BLLN makes it possible to view the Kaplan–Meier estimator as an asymptotic instance of the Bayesian MAP, and provides a probabilistic underpinning. The only available Bayesian view of the Kaplan–Meier estimator seems to be that of Susarla and van Ryzin [37]. In [37], a Dirichlet process prior was considered in a well-specified model, and it was shown there that the posterior mean converges to the Kaplan–Meier estimator as the parameter $\alpha$ of the Dirichlet process converges to 0.

**Acknowledgments.** We are grateful to Nicole Lazar who gave us valuable feedback and encouraged us to write this paper. Major impetus for improvement of previous version of this work came from reviewers, associate editor and main editor. Robert Niven [29] opened the possibility space to Section 2.4. We also want to thank Douglas Miller, Art Owen, Dale Poirier, Jing Qin and Giuseppe Ragusa for valuable comments and suggestions.

## REFERENCES

[1] Ben-Tal, A., Brown, D. E. and Smith, R. L. (1987). Posterior convergence under incomplete information. Technical Report 87-23, Univ. Michigan, Ann Arbor.

[2] Ben-Tal, A., Brown, D. E. and Smith, R. L. (1988). Relative entropy and the convergence of the posterior and empirical distributions under incomplete and conflicting information. Technical Report 88-12, Univ. Michigan.

[3] Brown, B. M. and Chen, S. X. (1998). Combined and least squares empirical likelihood. *Ann. Inst. Statist. Math.* **90** 443–450. MR1671990

[4] Cressie, N. and Read, T. (1984). Multinomial goodness-of-fit tests. *J. Roy. Statist. Soc. Ser. B* **46** 440–464. MR0790631

[5] Csiszár, I. (1998). The method of types. *IEEE Trans. Inform. Theory* **44** 2505–2523. MR1658767

[6] Csiszár, I. and Shields, P. (2004). Notes on information theory and statistics: A tutorial. *Found. Trends Comm. Inform. Theory* **1** 1–111.

[7] Dembo, A. and Zeitouni, O. (1998). *Large Deviations Techniques and Applications*, 2nd ed. Springer, New York. MR1619036

[8] Florens, J.-P. and Rolin, J.-M. (1994). Bayes, bootstrap, moments. Discussion Paper 94.13, Institute de Statistique, Université catholique de Louvain.

[9] Freedman, D. A. (1963). On the asymptotic behavior of Bayes' estimates in the discrete case. *Ann. Math. Statist.* **34** 1386–1403. MR0158483

[10] Freedman, D. A. (1999). On the Bernstein–von Mises theorem with infinite-dimensional parameters. *Ann. Statist.* **27** 1119–1140. MR1740119

[11] Ganesh, A. and O'Connell, N. (1999). An inverse of Sanov's Theorem. *Statist. Probab. Lett.* **42** 201–206. MR1680126

[12] Ghosal, A., Ghosh, J. K. and Ramamoorthi, R. V. (1999). Consistency issues in Bayesian nonanparametrics. In *Asymptotics, Nonparametrics and Time Series: A Tribute to Madan Lal Puri* 639–667. Dekker. MR1724692

[13] GHOSH, J. K. and RAMAMOORTHI, R. V. (2003). *Bayesian Nonparametrics*. Springer, New York. MR1992245

[14] GODAMBE, V. P. and KALE, B. K. (1991). Estimating functions: An overview. In *Estimating Functions* (V. P. Godambe, ed.) 3–20. Oxford Univ. Press, New York. MR1163993

[15] GRENDÁR. M. (2005). Conditioning by rare sources. *Acta Univ. M. Belii Ser. Math.* **12** 19–29. MR2258824

[16] GRENDÁR, M. and JUDGE, G. (2008). Large deviations theory and empirical estimator choice. *Econometric Rev.* **27** 513–525. MR2417666

[17] HJORT, N. L., MCKEAGUE, I. W. and VAN KEILEGOM, I. (2007). Extending the scope of empirical likelihood. *Ann. Statist.* To appear.

[18] IMBENS, G., SPADY, R. and JOHNSON, P. (1998). Information theoretic approaches to inference in moment condition models. *Econometrica* **66** 333–357. MR1612246

[19] JUDGE, G. G. and MITTELHAMMER, R. C. (2007). Estimation and inference in the case of competing sets of estimating equations. *J. Econometrics* **138** 513–531. MR2391321

[20] KEČKIĆ, J. D. and VASIĆ, P. M. (1971). Some inequalities for the gamma function. *Publ. Inst. Math. (Beograd) (N.S.)* **11** 107–114. MR0308446

[21] KERRIDGE, D. F. (1961). Inaccuracy and inference. *J. Roy. Statist. Soc. Ser. B* **23** 284–294. MR0123375

[22] KITAMURA, Y. and STUTZER, M. (1997). An information-theoretic alternative to generalized method of moments estimation. *Econometrica* **65** 861–874. MR1458431

[23] KITAMURA, Y. and STUTZER, M. (2002). Connections between entropic and linear projections in asset pricing estimation. *J. Econometrics* **107** 159–174. MR1889957

[24] KLEIJN, B. J. K. and VAN DER VAART, A. W. (2006). Misspecification in infinite-dimensional Bayesian statistics. *Ann. Statist.* **34** 837–877. MR2283395

[25] KULHAVÝ, R. (1996). *Recursive Nonlinear Estimation: A Geometric Approach. Lecture Notes in Control and Information Sciences* **216**. Springer, London. MR1402014

[26] LAZAR, N. (2003). Bayesian empirical likelihood. *Biometrika* **90** 319–326. MR1986649

[27] MITTELHAMMER, R., JUDGE, G. and MILLER, D. (2000). *Econometric Foundations*. Cambridge Univ. Press, Cambridge. MR1789434

[28] MONAHAN, J. F. and BOOS, D. D. (1992). Proper likelihoods for Bayesian analysis. *Biometrika* **79** 271–278. MR1185129

[29] NIVEN, R. K. (2005). Combinatorial information theory: I. Philosophical basis of cross-entropy and entropy. Available at arXiv:cond-mat/0512017.

[30] OWEN, A. (1988). Empirical likelihood ratio confidence interval for a single functional. *Biometrika* **75** 237–249. MR0946049

[31] OWEN, A. (2001). *Empirical Likelihood*. Chapman & Hall/CRC Press, New York.

[32] QIN, J. and LAWLESS, J. (1994). Empirical likelihood and general estimating equations. *Ann. Statist.* **22** 300–325. MR1272085

[33] RAGUSA, G. (2006). Bayesian likelihoods for moment condition models. Working paper, Univ. California, Irvine.

[34] SANOV, I. N. (1957). On the probability of large deviations of random variables. *Sb. Math.* **42** 11–44. (In Russian.) MR0088087

[35] SCHENNACH, S. (2005). Bayesian exponentially tilted empirical likelihood. *Biometrika* **92** 31–46. MR2158608

[36] SCHWARTZ, L. (1965). On Bayes procedures. *Z. Wahrsch. Verw. Gebiete* **4** 10–26. MR0184378

[37] SUSARLA, V. and VAN RYZIN, J. (1976). Nonparametric Bayesian estimation of survival curves from incomplete observations. *J. Amer. Statist. Assoc.* **71** 897–902. MR0436445
[38] WALKER, S. (2004). New approaches to Bayesian consistency. *Ann. Statist.* **32** 2028–2043. MR2102501
[39] WALKER, S. and DAMIEN, P. (2000). Practical Bayesian asymptotics. Working Paper 00-007, Business School, Univ. Michigan.
[40] WALKER, S., LIJOI, A. and PRÜNSTER, I. (2004). Contributions to the understanding of Bayesian consistency. Discussion Paper 13/2004, ICER, Applied Mathematics Series.
[41] WANG, D. and CHEN, S. X. (2007). Empirical likelihood for estimating equations with missing values. *Ann. Statist.* To appear.

DEPARTMENT OF MATHEMATICS
FACULTY OF NATURAL SCIENCES
BEL UNIVERSITY
TAJOVSKÉHO 40
SK-974 01 BANSKÁ BYSTRICA
SLOVAKIA
AND
INSTITUTE OF MATHEMATICS AND CS
SLOVAK ACADEMY OF SCIENCES (SAS)
INSTITUTE OF MEASUREMENT SCIENCES OF SAS
BRATISLAVA
SLOVAKIA
E-MAIL: marian.grendar@savba.sk

THE GRADUATE SCHOOL
UNIVERSITY OF CALIFORNIA
207 GIANNINI HALL
BERKELEY, CALIFORNIA 94720
USA
E-MAIL: judge@are.berkeley.edu